\documentclass[12pt]{amsart}
\usepackage{amscd,amssymb,verbatim}
\usepackage{amssymb,amsmath,amsthm,epsfig}

\newcommand{\N}{\mathbb{N}}

\newcommand{\R}{\mathbb{R}}
\newcommand{\vp}{\varepsilon}
\newcommand{\supp}{\text{supp}}
\newcommand{\nl}{\langle\hskip-1pt\langle}
\newcommand{\nr}{\rangle\hskip-1pt\rangle}
\newcommand{\nnl}{\langle\hskip-1pt\langle\hskip-1pt\langle}
\newcommand{\nnr}{\rangle\hskip-1pt\rangle\hskip-1pt\rangle}
\theoremstyle{plain}
\newtheorem{Thm}{Theorem}[section]
\newtheorem{Prop}[Thm]{Proposition}
\newtheorem{Lem}[Thm]{Lemma}

\newtheorem*{Rmk}{Remark}

\newtheorem*{Def}{Definition}

\begin{document}

\title{The Banach space \text{\boldmath$S$} is complementably minimal and
subsequentially prime}

\author{G. Androulakis and Th. Schlumprecht}

\thanks{Research of both authors was supported by NSF}

\subjclass{46B03, 46B20}


\maketitle

\noindent
{\bf Abstract} We first include a result of the second author showing
that the Banach space $S$ is complementably minimal. We then show that every 
block  sequence of the unit vector basis of $S$ has
a subsequence which spans a space  isomorphic to its square.
By the Pe{\l}czy\'nski decomposition method it follows that every
 basic sequence in $S$ which spans a space complemented 
in $S$ has a subsequence which spans a space isomorphic to $S$
(i.e. $S$ is a subsequentially prime space).

\bigskip

\section{Introduction}\label{sec1}

The Banach space $S$ was introduced by the second author as an example of an 
arbitrarily distortable Banach space \cite{S1}. In \cite{GM1} the
space $S$ was used to construct a Banach space which does not contain any
unconditional basic sequence. In this paper we are
concerned with the question  whether  or not $S$ is a prime space.
 We  will present two partial results: In Section 2 we show
that $S$ is complementably minimal, 
 and thereby answer a question of P.~G.~Casazza, 
 who asked whether or not $\ell_p$, $1\le p\le \infty$,  and $c_0$ are the only
 complementably minimal spaces. 
 In Section 3 we prove that $S$ is subsequentially prime.

Let us recall the above notions. A Banach space $X$ is called {\em
prime} \cite{LT} if every complemented infinite dimensional subspace
of $X$ is isomorphic to $X$. A. Pe{\l}czy\'nski \cite{P} showed that the
spaces $c_0$ and $\ell_p$ ($1 \leq p < \infty$) are prime, and
J. Lindenstrauss \cite{L} showed that this is also true for the space
$\ell_\infty$. New prime spaces were  constructed by W.T. Gowers and
B. Maurey \cite{GM2}. But it is still open whether or not $\ell_p$, $1\le p\le\infty$,  and $c_0$  
 are the only prime spaces with an unconditional basis.

A space $X$ is called {\em minimal} (a notion due
to H. Rosenthal) if every infinite dimensional subspace of $X$
contains a subspace isomorphic to $X$, and $X$ is called {\em
complementably minimal} \cite{CS} if every infinite dimensional
subspace of $X$ contains a subspace which is isomorphic to $X$ and
complemented in $X$. P.G. Casazza and E. Odell \cite{CO} showed that
Tsirelson's space $T$ \cite{T}, as described in \cite{FJ}, fails to
have a minimal subspace. On the other hand it was shown by
P.G. Casazza, W.B. Johnson and L. Tzafriri \cite{CJT} that the space
$T^*$ is minimal but not complementably minimal. Since $S$ is
complementably minimal, either $S$
is prime, or there exists a complemented subspace $X$ of $S$ such that
$X$ and $S$ give a negative solution to the Schroeder-Bernstein
problem for Banach spaces (see \cite{C} for a detailed discussion of
this question): if two space are isomorphic to complemented subspaces
of each other must they be isomorphic? Negative solutions to the
Schroeder-Bernstein problem for Banach spaces are given by W.T. Gowers
\cite{G}, and W.T. Gowers and B. Maurey  \cite{GM2}, but to our knowledge
 it is open whether  or not there are two Banach spaces $X$ and $Y$,
 both having an unconditional basis, so that $X$ is complemented in
 $Y$ and $Y$ is complemented in $X$,  but so that $X$ and $Y$ are not complemented.

  The following 
terminology was suggested to us by D.~Kutzarova.

\begin{Def}
A Banach space $X$
  with a basis is called {\em  subsequentially prime} 
if for every  basic 
sequence $(x_i)$ of $X$ such that the closed linear span of $(x_i)$ 
is complemented in $X$, there exists a subsequence $(y_i)$ such
that  the closed linear span of $(y_i)$ is isomorphic to $X$. 
\end{Def}

As mentioned above, we will show that the space $S$ is subsequentially prime.
 We do not know if $S$ is prime, we even do not know whether or not 
 the closed linear span of a block basis which is complemented in $S$ is
isomorphic to $S$.

We will need some notations. 
Let $c_{00}$ be the linear span of finitely supported real sequences, and
let $(e_i)$ denote its standard basis. 
 For $x\in c_{00}$,  $\supp(x)=\{i\in \N: x_i\not=0\}$ denotes the 
{\em support of}  $x$. For a finite  set $A$ the cardinality of $A$ 
 is denoted by $\# A$.
 If $E,F\subset \N$ we write $E<F$ if $\max E<\min F$, and
 we write $x<y$ for $x,y\in c_{00}$ if $\supp(x)<\supp(y)$. 
 A sequence $(x_i)_i$ in $c_{00}$ is a {\em block
sequence} of $(e_i)$ if $x_1<x_2<\ldots $.
 For $x=\sum_{i\in\N} x_i e_i\in c_{00}$ and $E\subset\N$ 
 $E(x)$ is the projection of $x$ onto the span of $(e_i)_{i\in E}$, i.e.
$E(x)=\sum_{i\in E} x_i e_i$. 
\smallskip

Recall \cite{S1} that the norm of $S$ is the unique norm on the completion of
$c_{00}$ which satisfies the implicit equation:
\begin{equation}\label{E:1.1}
\| x \| = \| x \|_{\ell_\infty} \vee \sup \limits _{ 
\stackrel{ \scriptstyle 
2 \leq n, \ E_i \subseteq \mathbb{N}, \ i=1, \ldots , n}
{\scriptstyle E_1< E_2 < \ldots <E_n}}
\frac{1}{f(n)} \sum_{i=1}^n \| E_i x \|
\end{equation}
where $\| \cdot \|_{\ell_\infty}$ denotes the norm of $\ell_\infty$ and
$f(n)=\log _2(n+1)$, for $n\in\N$.
For $x\in S$ and $\ell\in \mathbb{N}$, $\ell \ge2$, we define
$$
\| x\|_\ell := \sup_{E_1<E_2< \ldots < E_\ell} \frac{1}{f(\ell)}
\sum_{i=1}^\ell \| E_i(x)\|,
$$

We note that $\|\cdot\|_\ell$, $2\le \ell <\infty$, is an
equivalent norm on $S$ and we observe that for $x\in S$ and $2 \le \ell <
\infty$ we have
$$\frac{1}{f(\ell)} \|x\| \le\|x\|_\ell \le \|x\|\ \text{ and }\
\| x\| = \sup_{2\le \ell \le\infty} \|x\|_\ell\ .$$

Finally we put for any $2\le r<\infty$ and $x\in S$
$$||| x |||_r := \sup_ {\stackrel{\scriptstyle \ell\ge r}
{\scriptstyle \ell \in \mathbb{N} \cup \{ \infty \} } }\| x \|_\ell .$$

  Two sequences $(x_i)$,
$(y_i)$ in $S$ are called $c$-{\em equivalent},  for some $c \geq 1$,
 and we write $(x_i)\approx_c(y_i)$,
 if $\| \sum a_i x_i \| \stackrel{c}{\approx} \| \sum a_i y_i
\|$ for all $(a_i) \in c_{00}$, where for $c \geq 1$ and $a,b \geq 0$
we write $a \stackrel{c}{\approx} b$ to denote that $(1/c) a \leq b
\leq ca$.  If $(x_i)$ and $(y_i)$ are $c$-{\em equivalent} for some
 $c>1$
 we write $(x_i) \approx (y_i)$.   A basic sequence $(x_n)$ is called
 {\em $c$-subsymmetric}if it is $c$-unconditional and $c$-equivalent to
all of its subsequences.
 For two Banach spaces $X$ and $Y$  we write $X\approx_c Y$ if 
 there is an  isomorphism $T$ between $X$ and $Y$ with $\|T\|\cdot \|T^{-1}\|\le c$ 
and we write  $X\approx Y$ if  $X\approx_c$ for some $c\ge 1$.

If $(x_n)$ is a sequence in a Banach space $[x_n:n\in\N]$ denotes the closed linear
 span of $(x_n)$.
If not said otherwise, all statements in the following sections 
 refer to the space $S$.

\medskip
\noindent
 We would like to thank P.G. Casazza and D. Kutzarova
for valuable discussions.

\section{The Banach space $S$ is complementably minimal} \label{sec2}

The goal of this section is the proof of the 
 following Theorem.
\begin{Thm}\label{T:2.1} $S$ is complementably minimal.
\end{Thm}
 
 First
recall the following result which follows  from Lemma 5 of
\cite{S1}. 

\begin{Prop}\label{P:2.1} 
$\ell_1$ is block finitely represented in each block basis
of $(e_i)$, {\it i.e.}, if $\varepsilon >0$ and $m\in\mathbb{N}$, and if
$(y_n)$ is
a block basis of $(e_i)$ 
then there is a block basis $(z_i)_{i=1}^m$ of $(y_n)$ which is
$(1+\varepsilon)$-equivalent to the unit basis of $\ell_1^m$ (i.e.,
$\|\sum_{i=1}^m\alpha_i z_i\| \ge \frac{1}{1+\varepsilon} \sum_{i=1}^m
|\alpha_i|$ for $(\alpha_i)_{i=1}^m \subset \mathbb{R}$) 
\end{Prop}

The proof of the following  statement can be compiled from the proof
of Theorem~3 of \cite{S1}. Since the statement is crucial for 
our main result we include its proof.

\begin{Lem} \label{Lem:duo}
Let $\varepsilon >0$ and $\ell\in\mathbb{N}$. Then there is an $n=
n(\varepsilon,\ell)
\in\mathbb{N}$ with the following property: If $m\ge n$ and if $y= \frac{1}{m}
\sum_{i=1}^m x_i$, where $(x_i)_{i=1}^m$ is a normalized block basis
of $(e_i)$ which is $(1+\varepsilon/2)$-equivalent to the unit basis of
$\ell_1^m$, then
$$\sup_{E_1<E_2< \ldots < E_\ell} \sum_{i=1}^\ell \| E_i(y)\|
\le \| y\| +\varepsilon \le 1+\varepsilon\ .$$
\end{Lem}

\begin{proof}
Let $n\in\mathbb{N}$ so that $\frac{4\ell}{n} \le \varepsilon$ and assume $m\ge
n$ and $(x_i)_{i=1}^m$ are given as in the statement.
Furthermore, let $E_1< E_2 < \ldots < E_\ell$ be finite subsets of
$\mathbb{N}$.  Since $(e_i)$ is 1-unconditional we can assume that the
$E_j$'s  are
intervals in $\mathbb{N}$. This implies that for each $j\in
\{1,2,\ldots,\ell\}$
there are at most two elements $i_1,i_2\in \{1,\ldots,m\}$ so that
$E_j\cap \text{supp} (x_{i_s})\ne\emptyset$ and 
$\text{supp} \,(x_{i_s})\setminus E_j \ne \emptyset$, $s=1,2$.
For $j=1,2,\ldots,\ell$, let
$$ \tilde E_j = \cup \{ \text{supp} (x_i) : i\le m  \text{ and }
\text{supp} (x_i) \subset E_j \} .$$
It follows that if $y= \frac{1}{m} \sum_{i=1}^m x_i$, then $\| E_i(y)
- \tilde E_i (y)\| \le \frac{2}{m}$ and from the assumption that
$(x_i)_{i=1}^m$ is $(1+\varepsilon/2)$-equivalent to the
$\ell_1^m$-unit-basis we deduce that
\begin{align*}
\sum_{j=1}^\ell \| E_j(y)\|
& \le \frac{2\ell}{m} + \sum_{j=1}^\ell \| \widetilde E_j(y)\| 
= \frac{2\ell}{m} + \frac{1}{m} \sum_{j=1}^\ell \Big\| \sum_{\text{supp} (x_i)
\subset\widetilde E_j} x_i\Big\|\\
&\le \frac{2\ell}{m} + \frac{1}{m}\Big\| \sum_{i=1}^m x_i\Big\| 
(1+\varepsilon/2)
\le \varepsilon +\|y\|
\le \varepsilon +1\ 
\end{align*}
which finishes the proof.
\end{proof}

The  following theorem essentially proves that $S$ is minimal. We
postpone its proof.

\begin{Thm} \label{Thm:tria}
Let $(\varepsilon_n)\subset \mathbb{R}^+$ with $\sum \varepsilon_n
<\infty$  and
let $(y_n)$ be a normalized block basis of $(e_n)$ with the following
properties: There is a sequence $k_n\uparrow \infty$ in $\mathbb{N}$ so
that for all $n\in\mathbb{N}$
\begin{equation} \label{E:-ena}
\sup_{\stackrel{\scriptstyle k\le k_{n-1} }{\scriptstyle E_1<E_2<\ldots < E_k}}
\sum_{i=1}^k \| E_i(y_n)\| \le 1+\varepsilon_n
\end{equation}
\begin{equation} \label{E:-duo}
 \max\text{\em supp} (y_n) \le \varepsilon_n  f\Bigl( \frac{k_n}{3}\Bigr).
\end{equation}
Then $(y_n)$ is equivalent to $(e_n)$.
\end{Thm}

\noindent{\em Proof of Theorem \ref{T:2.1}}
 By  the usual
perturbation argument we only have to show that every block basis $(z_n)$
of $(e_n)$ has a further block basis which is equivalent to $(e_n)$.
Letting for example $\varepsilon_i = 2^{-i}$, $i=1,2,\ldots$, we have to
find
a normalized block $(y_n)$ of $(z_n)$ and a sequence $(k_n)$ in 
$\mathbb{N}$ so
that (\ref{E:-ena}) and (\ref{E:-duo}) of Theorem \ref{Thm:tria} are
satisfied. Indeed, put $k_0=1$
and  assume that
$k_0 < k_1 <\ldots < k_n$ and $y_1<y_2 <\ldots < y_n$ are already
defined for some $n\ge0$.
By Remark \ref{P:2.1}  and Lemma \ref{Lem:duo}  we can choose $y_{n+1}
>y_n$  in the linear
span of $(z_i)$ so that condition (\ref{E:-ena}) of Theorem \ref{Thm:tria} is
satisfied.   Since
$\lim_{i\to\infty} f(i)=\infty$ we then can choose $k_{n+1}$ so that (\ref{E:-duo})
is true.

In order to show that $S$ is complementably minimal we first observe that (\ref{E:1.1}) 
implies that 
every normalized block basis $(y_n)$ of $(e_n)$ dominates $(e_n)$, {\it i.e.},
that
 $\| \sum_{i=1}^\infty \alpha_i y_i \| \ge \| \sum_{i=1}^\infty
\alpha_i e_i\|$,
for all $(\alpha_i) \in c_{00}$. 
Secondly we apply the following  more general Proposition.\hfill$\square$

\begin{Prop} \label{Prop:tessera}
Let $Z$ be a Banach space with a $c_u$-unconditional basis
$(e_n)$, $c_u\ge 1$. Assume furthermore that there is a $c_d>0$ so that every
normalized block basis $(y_n)$ of $(e_n)$ $c_d$-dominates $(e_n)$
(i.e., $ c_d | \sum \alpha_i y_i\| \ge  \|\sum \alpha_i e_i\| $
for all $(\alpha_i)\in c_{00}$). Then a subspace of $Z$ generated by a
normalized block of $(e_n)$ which is equivalent to $(e_n)$ is
complemented in $Z$.
\end{Prop}

\begin{proof} W.l.o.g. we can assume that $(e_n)$ is a normalized and bimonotone
 basis of $Z$ (i.e. $\|[m,n](z)\|\le \|z\|$ for all $z=\sum_{i=1}^\infty z_ie_i\in Z$,
and $1\le m\le n$ in $\N$). Assume that $(y_n)$ is a block of $(e_i)$  which 
is $c_e$-equivalent  to $(e_i)$. 
Using the  assumption that $(e_n)$ is normalized and bimonotone we find $y^*_n\in Z^*$, for $n\in\N$,
 with $1=\|y^*_n\|=y^*_n(y_n)$ and $\supp(y^*_n)\subset[1+\max(\supp(y_{n-1})),\max(\supp(y_{n}))]$
 (where $y_0=0$ and $\max(\emptyset)=0$).
Define $T=\sum y^*_n\otimes y_n$, $x\mapsto \sum  y_n y^*_n(x)$ . We have to show that
$T$ is welldefined and bounded on $Z$, then it easily follows that it is a projection
 on $[y_n:n\in\N]$.
Let $x=\sum a_i e_i$, with $(a_i)\in c_{00}$. We can write 
 $x=\sum x_i=\sum \|x_i\| u_i,$   
with $x_i=[1+\max(\supp(y_{i-1})),\max(\supp(y_{i}))](x)$ and $u_i=x_i/\|x_i\|$ if
 $x_i\not=0$, and $u_i=e_{\max\supp y_i}$, otherwise.
 Then it follows that 
\begin{align*}
T(x)&=\Big\|\sum y_n y^*_n(x_n)\Big\|\le c_u\Big\|\sum \|x_n\| y_n \Big\|
 \le c_u c_e \Big\|\sum\|x_n\| e_n \Big\|\\
    & \le c_u c_e c_d \Big\|\sum \|x_n\|u_n \Big\|=c_u c_e c_d\|x\|,
\end{align*}
which proves the claim.
\end{proof}

 Let $x\in S$. If $\ell$ is the smallest element of $\N$, so that
 $\|x\|=\|x\|_{\ell}$  we call $\ell$ {\em the character of $x$} and
 write char$(x)=\ell$. If $\|x\|=\|x\|_{\ell_\infty}$ we write 
   char$(x)=\infty$.

The next Lemma makes  the following qualitative statement precise:
If $x\in c_{00}$, if $r>1$ is  ``big enough'', and if $E_1< E_2 < \ldots
< E_\ell$, $\ell \ge r$, are subsets of $\mathbb{N}$ so that
$$||| x |||_r = \frac{1}{f(\ell)} \sum_{i=1}^\ell \| E_i (x)\|\ ,$$
then for ``most of the $E_i$'s'' the character of $E_i(x)$ is ``much
bigger than $r$.''

\begin{Lem} \label{Lem:pente}
There is a constant $d>1$ so that for all $r\in \mathbb{R}_+$  with $f(r) >
d^2$,
$$||| x |||_r \le \left[ \frac{1}{1- \frac{d}{\sqrt{f(r)}}}\right]
\sup_{\stackrel{\scriptstyle \ell\ge r}{\scriptstyle E_1<E_2<\ldots < E_\ell}}
\frac{1}{f(\ell)} \sum_{i=1}^\ell ||| E_i(x) |||_{r^{f(r)}}$$
if $x\in c_{00}$ with $||| x |||_r \ne \|x\|_{\ell_\infty}$.
\end{Lem}

\begin{proof}
From the logarithmic behavior of $f$ we deduce that there is a constant $c>2$
so that the following inequalities hold
\begin{equation}\label{E1}
 f(\xi)-1\ge f(\xi)/c,\text{ whenever }\xi\ge 2,
\end{equation}
\begin{equation}\label{E2}
c f(\xi)\ge f(\xi\xi')-f(\xi),\text{ whenever }\xi,\xi'\ge c,
\end{equation}
\begin{equation}\label{E3}
f(\xi^{1/\sqrt{f(\xi)}})\le c \sqrt{f(\xi)},\text{ whenever }\xi\ge c,\text{ and }
\end{equation}
\begin{equation}\label{E4}
f(\xi^\nu)\le c \nu f(\xi) \text{ whenever }\xi\ge c\text{ and }\nu\ge1.
\end{equation}

 Choose $d= 4c^3$, let $r\in \mathbb{R}_+$ such that $f(r) > d^2$. In order
 to verify that this choice works 
 let  $x\in S$  with $|||x|||_r\not=\|x\|_{\infty}$. Let $\ell\ge r$
 and $E_1<E_2<\ldots<E_\ell$ so that
$$||| x |||_r = \frac{1}{f(\ell)} \sum_{i=1}^\ell \| E_i (x)\|\ .$$

For  $\tilde r, \widetilde R\in\R$, with 
$2\le \tilde r<\widetilde R$, let
 $M = M(\tilde r,\widetilde R) := \bigl\{ i\le \ell :\text{char}\, 
\bigl( E_i (x)\bigr) \in [\tilde r,\widetilde R[\ \bigr\}$
 and for $i\in M$ let $\ell_i \in [\tilde r,\widetilde R[$ be the
character of $E_i(x)$. We choose for each $i\in M$ finite subsets
of $E_i$, $E_1^i < E_2^i < \ldots < E_{\ell_i}^i $ so that
$$\| E_i(x)\| = \frac{1}{f(\ell_i)} \sum_{j=1}^{\ell_i} \| E_j^i (x)\|\ .$$
Now we observe that the set $\{ E_i: i\notin M\}\cup \bigcup_{i\in M}
\{ E_j^i :1\le j\le \ell_i\}$ is well ordered by ``$<$'' and its
cardinality is $\ell -\#\,M + \sum_{i\in M}\ell_i$ which is at least
$\ell$ and at most $\ell\tilde R$. Thus we deduce:
\begin{align}\label{E4a}
||| x |||_r
& = \frac{1}{f(\ell)} \sum_{i=1}^\ell \| E_i (x)\|\\
&\ge \frac{1}{f(\ell-\#\,M+ \sum\limits_{i\in M} \ell_i)}
\left[ \sum_{i=1,\, i\notin M}^\ell \| E_i(x)\| +
\sum_{i\in M}\ \sum_{j=1}^{\ell_i} \| E_j^i (x)\|\right]\notag\\
&\ge \frac{1}{f(\ell\widetilde R)} \left[ \sum_{i=1,\, i\notin M}^\ell
\| E_i(x)\| + \sum_{i\in M} f(\ell_i) \| E_i (x)\|_{\ell_i}\right]\notag\\
&\ge \frac{1}{f(\ell \widetilde R)} \left[ \sum_{i=1}^\ell
\| E_i(x)\| + \sum_{i\in M} \bigl( f(\tilde r)-1\bigr)
\| E_i (x)\| \right]  \notag  \\
&\ge \frac{1}{f(\ell\widetilde R)} \left[ \sum_{i=1}^\ell
\| E_i(x)\| + \frac{1}{c} f(\tilde r) \sum_{i\in M} \|E_i(x)\|\right]\  
   \text{ (using (\ref{E1}))}.\notag
\end{align}
Solving for $\frac{1}{f(\ell)} \sum_{i\in M} \|E_i(x)\|$ leads to the
following inequalities
\begin{align}\label{E:5}
\frac{1}{f(\ell)} \sum_{i\in M} \| E_i(x)\|
&\le  \frac{1}{f(\ell)} \left[ \frac{1}{f(\ell)} - \frac{1}{f(\ell
\widetilde R)}\right]  \frac{c f(\ell \widetilde R) }{
f(\tilde r)} \sum_{i=1}^\ell \| E_i (x)\|  \\
& = \frac{c}{f(\tilde r)} \frac{f(\ell\widetilde R) - f(\ell) }{
f(\ell)}\ ||| x |||_r\notag\\
&\le c^2 \frac{f(\widetilde R)}{f(\tilde r)f(\ell)} \
||| x |||_r
 \le c^2 \frac{f(\widetilde R)}{f(\tilde r) f(r)}\
||| x |||_r\  \text{ (using \ref{E2})}.\notag
\end{align}

Choosing for the pair of numbers $(\tilde r,\widetilde R)$ the values
$(2,r^{1/\sqrt{f(r)}})$, $(r^{1/\sqrt{f(r)}},r)$, $(r,r^{\sqrt{f(r)}})$,
and $(r^{\sqrt{f(r)}}, r^{f(r)})$ we deduce from the inequalities
 (\ref{E3}) and  (\ref{E4}) in each case that
 $\frac{f(\widetilde R)}{f(\tilde r) f(r)} \le \frac{c }{
\sqrt{f(r)}}$,
which implies together with (\ref{E:5}) that
$$\frac{1}{f(\ell)} \sum_{2\le \text{char}\, (E_i(x))< r^{f(r)}}
\|E_i(x)\| \le \frac{4c^3}{\sqrt{f(r)}} \ ||| x |||_r\ ,$$
and, thus, that
$$ ||| x |||_r \le \frac{d}{\sqrt{f(r)}} \ ||| x |||_r
+ \frac{1}{f(\ell)} \sum_{i=1}^\ell ||| E_i(x) |||_{r^{f(r)}},$$
yielding the lemma.
\end{proof}

\begin{Rmk} Note that in the proof of Lemma \ref{Lem:pente} the only properties
 of the function $f$ which  was needed was
that fact that it was increasing and   that there is a $c>2$ so that
 the inequalities (\ref{E1}) - (\ref{E4}) hold.  

Thus if $\ell_0\in\N$ and $g:[\ell_0,\infty)\to (1,\infty)$ is an increasing function
 so that there is a $c>\ell_0$ for which (\ref{E1}) (whenever $\xi\ge \ell_0$),
(\ref{E2}),(\ref{E3}) and (\ref{E4}) hold then the conclusion of Lemma \ref{Lem:pente} holds
for the completion of  $c_{00}$ under the norm $\nl\cdot\nr$ defined implicitly by
$$\nl x\nr=\|x\|_{\ell_\infty}\vee\sup_{\ell\ge\ell_0,E_1<E_2<\ldots E_\ell}\frac1{g(\ell)}
 \sum_{i=1}^\ell \nl E_i(x)\nr, \text{ whenever }x\in c_{00}.$$  
\end{Rmk}

\begin{proof}[Proof of Theorem \ref{Thm:tria}]
Let $(y_n)$ , $(k_n)$  and $(\varepsilon_n)$ be given as in the statement
of Theorem \ref{Thm:tria} and let $d\ge1$ be as in Lemma
\ref{Lem:pente}.  
For $r\ge1$  we put $r_0:= r$ and, assuming $r_k$ was already defined,
we let $r_{k+1} = r_k^{f(r_k)}$. From the properties of the function $f$
it follows that there is an $R>1$ so that the value
\begin{equation} \label{eksi}
\beta (r) := \prod_{k=0}^\infty \left( \frac{1}{1- \frac{d}{\sqrt{f(r_k)}}}
\right)  \frac{f(9 r_k) }{f(r_k)}
\end{equation}
is finite whenever $r\ge R$.
By induction we will show that for every $m\in \mathbb{N}$  
 and every $(\alpha_i)_{i=1}^m\subset \R$.
\begin{equation} \label{epta}
||| \sum_{i=1}^m \alpha_i y_i |||_{\,r} \le \beta (r) \max_{i_0\ge 1}
  \left[ |\alpha_{i_0} | + \Big\|\sum_{i>i_0}
\alpha_i e_i\Big\|+ \sum_{i=1}^m |\alpha_i|\varepsilon_i\right]\ ,
\end{equation}

Since $||| \cdot |||_r $ is equivalent to $\|\cdot\|$ for all
$r\ge1$, since $\sum \varepsilon_i <\infty$, and since $\| \sum\alpha_i
y_i\|
\ge\|\sum \alpha_ie_i\|\ge \max_{i\in\mathbb{N}} |\alpha_i|$ for $(\alpha_i)
\in c_{00}$ this would prove the assertion of Theorem \ref{Thm:tria}.

For $m=1$ the claim is trivial. Assume it is true for all
positive integers smaller
than some $m>1$ and let $r\ge R$, $(\alpha_i)_{i=1}^m\in c_{00}$.
 Let $y= \sum_{i=1}^m \alpha_iy_i$.
If $||| y |||_r = \|y\|_{\ell_\infty}$ the assertion follows easily
since $\|y\|_{\ell_\infty} \le \max_{i\le m}|\alpha_i|$. Otherwise we can use
Lemma \ref{Lem:pente} in order to find an $\ell\ge r$ and finite
subsets of  $\mathbb{N}$,
$E_1<E_2<\ldots < E_\ell$ so that (with $\gamma (r) = 1/(1- {d}/{\sqrt{f(r)}})$)
\begin{equation} \label{E:okto}
||| y |||_r \le \gamma (r)  \frac{1}{f(\ell)}
\sum_{j=1}^\ell ||| E_j(y) |||_{r^{f(r)}} .
\end{equation}

We can assume that for all $j\le \ell$, $E_j\subset\bigcup_{i=1}^m \supp(y_i)$. 
 For $j=1,2,
\ldots,\ell$ we put $E_j^1 := E_j \cap \text{supp} (y_{s(j)})$, $E_j^2 :=
E_j \cap \text{supp} (y_{t(j)})$ and $E_j^3 = E_j\backslash (E_j^1 \cup
E_j^2)$ where $s (j) := \min \{ i: E_j\cap \text{supp} (y_i) \ne
\emptyset\}$
and $t(j) : = \max \{ i:E_j\cap \text{supp} (y_i) \ne \emptyset\}$. We put
$\widetilde {\mathcal E} := \{ E_j^1,E_j^2,E_j^3$, $j\le \ell\}
\backslash \{\emptyset\}$
and note that $\widetilde {\mathcal E}$ can be ordered into 
$\widetilde {\mathcal E} = \{
\widetilde E_1,\widetilde E_2,\ldots,\widetilde E_{\tilde\ell}\}$ with
$\ell \le \tilde\ell \le 3\ell$ and $\widetilde E_1 <
\widetilde E_2 <\ldots < \widetilde E_{\tilde \ell}$.

Secondly we observe that $\widetilde {\mathcal E}$ can be partitioned
into  $m+1$ sets
$\widetilde {\mathcal E}_0,\widetilde {\mathcal E}_1,\ldots,
\widetilde {\mathcal E}_m$ defined in the following way:
$\widetilde {\mathcal E}_0 := \{ E\in \widetilde {\mathcal E} : E$ fits
with  $(\text{supp} (y_i))_{i\in\mathbb{N}}\}$
(where we say that $E$ fits with a sequence $(A_n)$ of disjoint subsets
of $\mathbb{N}$ if and only if for all $n$, $E\cap A_n\ne\emptyset$
implies  that $A_n\subset E$) and for $1\le i\le m$ we let
$\widetilde{\mathcal E}_i := \bigl\{ E: E\in \widetilde {\mathcal E} \
\text{ and }\
E\subsetneqq \text{supp} (y_i)\big\}$.

For $i=1,\ldots,m$ we let $\ell_i := \#\,\widetilde {\mathcal E}_i$
(note that  $\widetilde {\mathcal E}_i$
may be empty) and let $i_0=1$ if for all $i\le m$, $\ell_i\le k_{i-1}$
otherwise put $i_0 := \max \{ i\le m:\ell_i >k_{i-1}\}$.

From (\ref{E:okto}) we deduce now that (recall that $r_1=r^{f(r)}$)
\begin{align} \label{E:ennia}
\qquad ||| y |||_r
& \le \frac{\gamma (r)}{f(\ell)} \sum_{j=1}^{\tilde \ell} |||
\widetilde E_j (y) |||_{r_1}\\
&\le \frac{\gamma (r) }{f(\ell)} \Biggl[ \sum_{j=1}^{\tilde\ell}
||| \widetilde E_j \biggl( \sum_{i<i_0} \alpha_iy_i\biggr) |||_{r_1}
+ \sum_{j=1}^{\tilde\ell} ||| \widetilde E_j (\alpha_{i_0} y_{i_0})
|||_{r_1} \nonumber \\
&\qquad\qquad + \sum_{j=1}^{\tilde\ell} ||| \widetilde E_j \biggl( \sum_{i>i_0}
\alpha_i y_i\biggr) |||_{\,r_1} \Biggr] \nonumber \\
&\le \frac{\gamma (r)}{f(\ell)} \sum_{i<i_0}\ \sum_{j=1}^{\tilde\ell}
|\alpha_i| \cdot \| E_j(y_i)\|
+   \frac{\gamma (r) f(\tilde \ell)}{f(\ell)}
  \frac{|\alpha_{i_0}|}{f(\tilde\ell)} \sum_{j=1}^{\tilde\ell}
\| \widetilde E_j (y_{i_0})\|  \nonumber\\      
&\qquad + \frac{\gamma (r)}{f(\ell)} \sum_{i>i_0,\, {\mathcal E}_i \ne\emptyset}
|\alpha_i| \sum_{E\in {\mathcal E}_i} \|E(y_i)\|
+ \frac{\gamma (r)}{f(\ell)} \sum_{E\in {\mathcal E}_0}
||| E\biggl( \sum_{i>i_0} \alpha_i y_i\biggr) |||_{r_1}. \nonumber
\end{align}
If $i_0\ne 1$ we deduce that the first term in the above sum,
can be estimated as follows (we use
condition (\ref{E:-duo}) of the statement of Theorem \ref{Thm:tria} and
note that  from the choice of
$i_0$ it follows that $\ell\ge\tilde\ell /3\ge \ell_{i_0} /3 \ge k_{i_0-1} /3
\ge k_i/3$ for $i<i_0$):
\begin{equation*}
\frac{\gamma (r)}{f(\ell)} \sum_{i<i_0} \
\sum_{j=1}^{\tilde\ell} |\alpha_i| \cdot \| E_j(y_i)\|
\le \gamma (r) \sum_{i<i_0} \frac{|\alpha_i| }{f(k_i/3)} \cdot
\#\, \text{supp} (y_i)
\le \gamma (r) \sum_{i=1}^{i_0-1} \varepsilon_i |\alpha_i|\ .
\end{equation*}
The second term can be estimated as follows:
\begin{equation*}
\frac{\gamma (r) f(\tilde\ell)}{f(\ell)}
 |\alpha_{i_0}| \frac{1}{f(\tilde \ell)}
\sum_{j=1}^{\tilde \ell} \|\widetilde E_j (y_{i_0})\|
\le \frac{\gamma (r) f(3\ell)}{f(\ell)} |\alpha_{i_0}|
\, \| y_{i_0}\| 
\le \frac{\gamma (r) f(3r)}{f(r)}  |\alpha_{i_0}|.
\end{equation*}

By condition (\ref{E:-ena}) of the statement of Theorem \ref{Thm:tria}
and the  definition
of $i_0$ we deduce that
 $\sum_{E\in {\mathcal E}_i}\| E(y_i)\| \le 1+\varepsilon_i \ ,\ \text{ if }
\ i  >i_0\ \text{ and }\ {\mathcal E}_i \ne \emptyset\ .$
Thus, we observe for the third term that
 $$\frac{\gamma (r)}{f(\ell)} \sum_{i>i_0,\, {\mathcal E}_i \ne\emptyset}
|\alpha_i| \sum_{E\in {\mathcal E}_i} \|E(y_i)\|
\le\frac{\gamma (r)}{f(\ell)}  \sum_{i>i_0,\, {\mathcal E}_i\ne\emptyset}
 (1+\varepsilon_i)  |\alpha_i|.$$
For the last term we  apply the induction hypothesis and find for each
$E\in {\mathcal E}_0$ an $i_E \in \{ i>i_0 :\text{supp} (y_i) \subset
E\}  \cup \{0\}$
so that
$$\sum_{E\in {\mathcal E}_0} ||| E\biggl( \sum_{i>i_0}
\alpha_i y_i\biggr) |||_{\,r_1}
\le \beta (r_1)  \sum_{E\in{\mathcal E}_0}
\Biggl[ |\alpha_{i_E} | + \Big\|\sum_{\stackrel{\scriptstyle i>i_E}{
\scriptstyle \text{supp} (y_i) \subset E} } \alpha_i e_i\Big\|
+ \sum_{\text{supp} (y_i)\subset E} |\alpha_i|\varepsilon_i\Biggr]\ .$$
Let 
$
{\mathcal A}= \bigl\{ \{i\} : i\!>\!i_0, \, {\mathcal E}_i\ne\emptyset\bigr\}
\cup \bigl\{ \{i_E\}: E\in {\mathcal E}_0 \bigr\}
\cup \bigl\{ \{ i\!>\!i_E : \text{supp} (y_i)\!\subseteq \!E\} : E\!\in\!
{\mathcal E}_0\bigr\} \backslash \{\emptyset\}
$
and note that ${\mathcal A}$ consists of subsets of $\{i_0+1,i_0+2\ldots\}$ 
has at most $3\tilde \ell \le 9\ell$ elements and is
well ordered by $<$. Finally we deduce from (\ref{E:ennia}) and the
above  estimates that
\begin{align*}
 ||| y|||_r
&\le \gamma (r) \sum_{i<i_0} |\alpha_i| \varepsilon_i +
\frac{\gamma (r) f(3r)}{f(r)} |\alpha_{i_0}| 
+ \frac{\gamma (r)}{f(\ell)} \sum_{i>i_0,\, {\mathcal E}_i\ne
\emptyset}^m
|\alpha_i| + \frac{\gamma (r)}{f(\ell)}
\sum_{i=1,\, {\mathcal E}_i\ne\emptyset}^m \varepsilon_i |\alpha_i|\\
&\qquad + \frac{\gamma (r)\beta (r_1)}{f(\ell)} \sum_{E\in{\mathcal E}_0}
|\alpha_{i_E}| 
+ \frac{\gamma (r)\beta (r_1)}{f(\ell)} \sum_{E\in {\mathcal E}_0}
\Big\| \sum_{i>i_E,\, \text{supp} (y_i)\subset E} \alpha_i e_i\Big\|\\
&\qquad + \frac{\gamma (r)\beta (r_1)}{f(\ell)} \sum_{\stackrel{\scriptstyle
\text{supp} (y_i)\subset E}{\scriptstyle E\in {\mathcal E}_0} }|\alpha_i|
\varepsilon_i\\ &\le \beta (r) \left[ |\alpha_{i_0}| + \frac{f(r)}{f(\ell)f(9r)}
\sum_{A\in {\mathcal A}} \Big\| A\biggl( \sum_{i=1}^m \alpha_ie_i
\biggr)\Big\|
+ \sum_{i=1}^m \varepsilon_i |\alpha_i| \right]\\
&\text{(Note that }\ \beta (r) = \beta (r_1) \gamma(r) 
f(9r)/f(r))\\
&\le \beta (r) \left[ |\alpha_{i_0}| + \frac{f(r)}{f(9r)} 
\frac{f(9\ell)}{f(\ell)} \Big\|\sum_{i=1}^m \alpha_i e_i\Big\|
+ \sum_{i=i_0+1}^m \varepsilon_i |\alpha_i|\right]
\text{ (since }\#\, {\mathcal A} \le 9\ell)\\
&\le \beta (r) \left[|\alpha_{i_0}| +\Big\| \sum_{i=1}^m \alpha_i e_i
\Big\| + \sum_{i=i_0+1}^m \varepsilon_i |\alpha_i|\right]
\text{ (since }\ell \ge r).
\end{align*}
This proves the induction step and completes the proof of Theorem
\ref{Thm:tria}.
\end{proof}

\section{The Banach space $S$ is subsequentially prime} \label{S:3}

The main result of this section is the following Theorem.
\begin{Thm}\label{T:3.1}
The space $S$ is subsequentially prime.
\end{Thm}
Theorem \ref{T:3.1} will essentially follow from 
A.~Pe{\l}czy\'nski's decomposition method and 
 the following theorem. 

\begin{Thm}\label{T:3.2}
Let $(x_i)$ be a normalized block sequence of $(e_i)$ in $S$
 and let $(k_i)$ be a subsequence of $\N$. There exists a
subsequence $(y_i)$ of $(x_i)$ with the following property:

If $(s_i)$ and $(t_i)$ are strictly increasing in $\N$, and $s_i,t_i\le k_i$, for $i\in\N$, then
 $(y_{s_i})$ and  $(y_{t_i})$ are  equivalent.
\end{Thm}

Before giving the proof of Theorem \ref{T:3.2} we  need a  result, for which we
introduce 
 the following norm $\nl\cdot\nr$ on $c_{00}$.

For $x\ge 3$ define $g=\log_2(1+\frac x2)$ and let $\nl\cdot\nr$ be the norm  which is 
 implicitly defined by
$$\nl x\nr=\|x\|_{\ell_\infty}\vee\sup_{\ell\ge 3,E_1<E_2<\ldots E_\ell}\frac1{g(\ell)}
 \sum_{i=1}^\ell \nl E_i(x)\nr, \text{ whenever }x\in c_{00}.$$  

\begin{Lem}\label{L:3.3}
 The norm $\nl\cdot \nr$ is  equivalent to the norm $\|\cdot\|$ on $S$.
\end{Lem}

\begin{proof} First note that since $f(\ell)\ge g(\ell+1)$ whenever $\ell\ge 2$ it
 follows that $\nl\cdot\nr\ge  \|\cdot\|$.  Thus, we only have to show that
for some constant $C$ it follows that  $\nl\cdot\nr\le C \|\cdot\|$. The proof
will be similar to the proof of Theorem \ref{Thm:tria}.

For $\ell$ and $r\in[3,\infty)$ let $\nl\cdot\nr_\ell$ and $\nnl\cdot\nnr_r$  be defined 
 as $\|\cdot\|_\ell$ and $|||\cdot|||_r$ respectively.
 Let $\tilde d>1$ be chosen so that the statement of Lemma \ref{Lem:pente} holds
 (see Remark after proof of  Lemma \ref{Lem:pente}),
i.e. so that 
\begin{equation}\label{E:3.3.1}
\nnl x \nnr_r \le \left[ \frac{1}{1- \frac{\tilde d}{\sqrt{g(r)}}}\right]
\sup_{\stackrel{\scriptstyle \ell\ge r}{\scriptstyle E_1<E_2<\ldots < E_\ell}}
\frac{1}{g(\ell)} \sum_{i=1}^\ell \nnl E_i(x) \nnr_{r^{g(r)}}
\end{equation}
if $r\in \mathbb{R}_+$  with $g(r) > \tilde d^2$, and
$x\in c_{00}$ with $\nnl x _r \nnr\ne \|x\|_{\ell_\infty}$.

For $r\in[3,\infty)$ we define $r_0=r$ and, assuming that $r_k$ had been defined, let
 $r_{k+1}=r_k^{g(r_k)}$.  Then we deduce that there is an $R\ge 3$ so that 
 for all $r\ge R$ 
$$\tilde\beta(r)=\prod_{k=0}^\infty \frac1{1-\frac{\tilde d}{\sqrt{g(r_k)}}}\frac{g(2r_k)}{g(r_k)}$$
is finite.
By induction on   $m\in\N$  we proof that for each $x\in c_{00}$ so that
 $\#\supp(x)\le m$  and $r\ge R$ it follows
that

\begin{equation}\label{E:3.3.2}
\nnl x \nnr_r\le\tilde\beta(r) \|x\|.
\end{equation}

Assume that (\ref{E:3.3.2}) is true for all $z\in c_{00}$ for which
 $\#\supp(z)< m$  and assume that $x \in c_{00}$ with $\#\supp(x)=m$.
If $\nnl x \nnr_r=\|x\|_{\ell_\infty}$ the claim follows immediately. Otherwise it follows
 from (\ref{E:3.3.1})  that  for some $\ell \ge r$, $\ell\in\N$, and
 some choice of $E_1<E_2<\ldots E_\ell$ we have 
\begin{align*}
\nnl x \nnr_r&\le \left[ \frac{1}{1- \frac{\tilde d}{\sqrt{g(r)}}}\right]
\frac{1}{g(\ell)} \sum_{i=1}^\ell \nnl E_i(x) \nnr_{r^{g(r)}}\\
&\le\tilde\beta(r^{g(r)})\left[ \frac{1}{1- \frac{\tilde d}{\sqrt{g(r)}}}\right]
 \frac{1}{g(\ell)}\sum_{i=1}^\ell \| E_i(x) \| \text{ (By  the induction hypothesis)}\\
&=\tilde\beta(r) 
\frac{g(r)}{g(2r)} \frac{1}{g(\ell)}\sum_{i=1}^\ell \| E_i(x) \| 
 \le\tilde\beta(r) \|x\| 
\text{ (Since $\frac{g(r)}{g(2r)}
         \le\frac{g(\ell)}{g(2\ell)}=\frac{g(\ell)}{f(\ell)}$)}
\end{align*}
which finishes the induction step and the proof of Lemma \ref{L:3.3}
\end{proof}

\begin{proof}[Proof of Theorem \ref{T:3.2}]

First we note that we can assume that $\lim_{n\to\infty }\|x_n\|_{\ell_\infty}=0$. Indeed,
  from Theorem \ref{Thm:tria} it follows that there is  a normalized
 block $(z_k)$ in $S$ which is equivalent to the unit vector basis $(e_k)$ and
has the property that $\lim_{n\to\infty }\|z_n\|_{\ell_\infty}=0$. Thus, we could replace 
 each $x_n$  by the vector $x_n'$  in the span of  $(z_k)$ whose
coordinates with respect to the $z_k$'s are the coordinates of $x_n$ with respect to the $e_k$'s.

 Let $\vp>0$, and $y\in c_{00}$ with $\|y\|=1$ and $\|y\|_\infty\le \vp/2$.
  We can write $y$ as 
 $y=\sum_{i=1}^{\ell(y,\vp)} y(i,\vp),$ 
where $\ell(y,\vp)\in\N$ and
$y(1,\vp)<y(2,\vp)<\ldots y(\ell(y,\vp),\vp)$ such that
 $\| y(i,\vp)\|\le \vp$ for $i=1,2,\ldots,\ell(y,\vp)$. 
 Further more we could choose the
$y(i,\vp)$'s to have maximal support in the following sense.
First we choose $y(i,\vp)=[1,n_1](y)$ with $n_1\in \N$  being the largest  $n\in\N$, $n\le \max\supp (y)$,
 so that $\|[1,n](y)\|\le \vp$, then we choose $n_2>n_1$
 being the largest  $n\in\N$, $n\le\supp(y)$, so that $\|[n_1+1,n_2](y)\|\le \vp$. 
 We can continue this way
until we have exhausted the support of $y$. This  defines $\ell(y,\vp)$ 
and the vectors $(y(i,\vp))_{i=1}^{\ell(y,\vp)}$ uniquely and from the assumption
 that $\|y\|_\infty\le \vp/2$ it follows that $\vp/2\le \|y(i,\vp)\|$ for all
 $i\in\{1,2,\ldots \ell(y,\vp)-1\}$. From the definition of the norm of $S$ is
follows  that
\begin{equation*}
\frac{\ell(y,\vp)-1}{f(\ell(y,\vp))}\cdot\frac{\vp}{2}\le
\frac1{f(\ell(y,\vp))}\sum_{i=1}^{\ell(y,\vp)}\|y(i,\vp)\|\le\|y\|=1\le
\sum_{i=1}^{\ell(y,\vp)}\|y(i,\vp)\|\le \ell(y,\vp)\vp.
\end{equation*}
Thus for $\vp>0$ there are two numbers $H(\vp)\ge h(\vp)$ in $\N$,
 with $h(\vp)\nearrow \infty$, if $\vp\searrow 0$, and so that for any
$y\in c_{00}$ with   with $\|y\|=1$ and $\|y\|_\infty\le \vp/2$ it follows that
  $h(\vp)\le \ell(y,\vp)\le H(\vp)$.
 
We  now apply this ``splitting procedure'' to the elements of our sequence
 $(x_n)$. By induction on $n\in\N$ we find an
 infinite  subset $M_n$ of $\N$, with
 $\N\supset M_1\supset M_2\supset\ldots M_n$ and 
 $\min M_1<\min M_2<\ldots \min  M_n$, numbers $\vp(n)>0$ and
  $p(n)\in\N$ satisfying the following  three properties.
\begin{align}
\label{E:3.2.1} &\text{For all $m\in M_n$ we have $\ell(\vp(n),x_m)=p(n)$}.\\
\label{E:3.2.2} &\text{For any choice of $n\le s_0<s_1<\ldots s_{k_n+1}$ and $n\le t_0<t_1<\ldots t_{k_n+1}$
  in $M_n$}\\
 &\text{ it follows that }
(x_{s_j}(i,\vp(n)))_{0\le j\le k_n+1 ,i\le p(n)}\approx_{1+\vp(n)}
(x_{t_j}(i,\vp(n)))_{0\le j\le k_n+1 ,i\le p(n)} .\notag\\
\label{E:3.2.3}&\text{If $n>1$ it follows that  
 $\sum_{1=i}^{ n-1}\|x_{\min M_i}\|_{\ell_1}< f(h(\vp(n)))2^{-n}$
 and $p(n-1)\vp(n) <2^{-n}$}.
\end{align}
For $n=1$ we simply choose $\vp(1)=1$ (thus $\ell(x_n,\vp(1))=1$), $p(1)=1$, and using
  compactness and the usual stabilization
 argument we can pass to a subsequence $M_1$ of $\N$ so that (\ref{E:3.2.2}) holds.
Assuming we made our choices  of $M_j$, $\vp(j)$, and $p(j)$, for all $j< n$, we
first choose $\vp(n)$ so that (\ref{E:3.2.3}) is satisfied 
 (recall that $h(\vp)\nearrow\infty$, for $\vp\searrow0$), and then again
  using compactness and the usual stabilization argument we
 can pass to an $M_n\subset M_{n-1}\setminus\{\min M_{n-1}\}$  and find  a $p(n)\in\N$ so that
  $\|x_m\|\le\vp(n)/2$ whenever $m\in M_n$ and so that
 (\ref{E:3.2.1}) and  (\ref{E:3.2.2}) are satisfied.

  For $n\in\N$ we now define $y_n=x_{\min(M_n)}$
  and  $y_n(i,j)=y_n(i,\vp(j))$ if $i\le n$ and $j\le p(i)$,
 and prove the following claim 
 by induction on $N\in\N$:

\noindent{\bf Claim:} For every $n\in\N$, every 
$K,L\in\{1,2,\ldots,p(n)\}$, every 
 $(\alpha_i)_{i=0}^{n+N+1}\subset\R$, and
 every choice of $n\le s_0 <s_1 \ldots s_{n+N+1}$ and $n\le t_0 <t_1 \ldots t_{n+N+1}$
with $s_i,t_i\le k_{i+n}$, for $i=0,1\ldots,n+N+1$,
 it follows for
\begin{align}\label{E:3.2.4}
 &x=\alpha_0\sum_{j=K}^{p(n)} y_{s_0}(n,j) +\sum_{i=1}^N \alpha_i y_{s_i} +
  \alpha_{N+n+1}\sum_{j=1}^{L} y_{s_{n+N+1}}(n,j), \text{ and }\\
\label{E:3.2.5} &\tilde x= \alpha_0\sum_{j=K}^{p(n)} y_{t_0}(n,j) +\sum_{i=1}^N \alpha_i y_{t_i} +
  \alpha_{N+n+1}\sum_{j=1}^{L} y_{t_{n+N+1}}(n,j)
\end{align}
that
\begin{equation}\label{E:3.2.5a}
\|x\|\le c(n)\bigl[\max_{E<F}\nl E(\tilde x)\nr+\nl F(\tilde x)\nr\bigr]
\end{equation}
where
\begin{equation}\label{E:3.2.5b}
c(n)=\sum_{i=n}^N 2^{-i}+\vp(i)
\end{equation}
Since by Lemma \ref{L:3.3}  $\nl \cdot\nr$  is an equivalent norm
 on $S$ (and therefore also the norm 
$c_{00}\ni x\mapsto\max_{E<F}\nl E(x)\nr+\nl F(x)\nr$ has this property) the
claim implies  the theorem.

For $N=0$ the claim follows directly from (\ref{E:3.2.2}). Assume the claim to be
true for $\tilde N<N$ and let $x$ and $\tilde x$ be given as in 
(\ref{E:3.2.4}) and (\ref{E:3.2.5}). 

We choose $\ell\in\N$ so that $\|x\|=\|x\|_\ell$.
If $\ell<p(k_n)$ we let $i_0=0$, and other wise we choose
$i_0\in\N\cup\{0\}$ so that $p(k_{i_0+n})\le \ell <p(k_{i_0+n+1})$.
We split $x$ into three vectors $x^{(1)}$, $x^{(2)}$ and $x^{(3)}$ as follows.
 If $i_0=0$ we let $x^{(1)}=0$, otherwise we put
 \begin{equation}\label{E:3.2.6}
x^{(1)}=\alpha_0\sum_{j=K}^{p(n)} y_{s_0}(n,j) +\sum_{i=1}^{(i_0-1)\wedge N} \alpha_i y_{s_i}.
\end{equation}
and $\tilde x^{(1)}$ is defined as $x^{(1)}$, by replacing the $s_i$'s by $t_i$'s.
 From (\ref{E:3.2.3})  and the  choice of $i_0$ it follows that
 (note that $h(\vp(s_{i_0-1}+1))\le h(\vp(s_{i_0}))\le p(s_{i_0})\le  p(k_{n+i_0})\le \ell$)
  
\begin{equation}\label{E:3.2.7}
\|x^{(1)}\|_\ell\le \max_{0\le i\le n+N}|\alpha_i|\frac1{f(\ell)}
   \sum_{i=0}^{(i_0-1)\wedge N} \|y_{s_i}\|_{\ell_1}
\le \|\tilde x\| \frac{f(h(\vp(s_{i_0})))}{f(\ell)}\le \|\tilde x\| 2^{-n-i_0}.
\end{equation}
Secondly, if we let
\begin{equation}\label{E:3.2.8}
x^{(2)}=\begin{cases}
\alpha_0\sum_{j=K}^{p(n)} y_{s_0}(n,j) +\sum_{i=1}^{(k_n+1)\wedge N} \alpha_i y_{s_i}
 &\text{if $i_0=0$}\\
{\phantom{--}}
 \sum_{i=(i_0-1)\wedge N+1}^{(k_{n+i_0}+1)\wedge N} \alpha_i y_{s_i} &\text{if $i_0\not=0$}
\end{cases}
\end{equation}
 and  $\tilde x^{(1)}$ is defined as $x^{(1)}$, by replacing the $s_i$'s by $t_i$'s.
 We deduce from (\ref{E:3.2.2}) that
\begin{equation}\label{E:3.2.9}
\|x^{(2)}\|\le(1+\vp(s_{i_0}))\|\tilde x^{(2)}\|
 \le (1+\vp({n+i_0}))\|\tilde x^{(2)} \|. 
\end{equation}
Finally we let 
\begin{equation}\label{E:3.2.10}
x^{(3)}=x-x^{(1)}-x^{(2)}=\sum_{ i=k_{n+i_0}+2}^N\alpha_i y_{s_i}+
 \alpha_{N+n+1}\sum_{j=L}^{p(n)} y_{s_{n+N+1}}(n,j)
\end{equation}
and, again, define $\tilde x^{(3)}$ as $x^{(3)}$, by replacing the $s_i$'s by $t_i$'s.
 Choose $E_1<\ldots E_\ell$ so that
\begin{equation*} 
\|x^{(3)}\|_\ell=\frac1{f(\ell)}\sum_{j=1}^\ell \|E_j(x^{(3)})\|.
\end{equation*}

Note that for any  $i\ge k_{n+i_0}+2$ it follows that
 $s_i\ge  k_{n+i_0}+2+n\ge k_{n+i_0}+2$
 and by  (\ref{E:3.2.3}) it follows that 
$\ell\vp(k_{n+i_0}+2)\le p(k_{n+i_0}+1)\vp(k_{n+i_0}+2)<2^{-k_{n+i_0}-1}\le 2^{-n-i_0-1}$.
Let $n'=k_{n+i_0}+2+n$. For each $j=1,2\ldots \ell$
 we can perturb $E_j$  into a set $F_j$ (by  possibly taking some part of the support of some $y_{s_i}(n',u)$
 away at the beginning of $E_j$ and adding some part of the support of some $y_{s_j}(n',v)$
   at the end of $E_j$)
 so that for each $j=1,2\ldots \ell$, each   $i\ge k_{n+i_0}+2$,
and each $u=1,2,\ldots p(n')$
 $F_j$ either contains $\supp (y_{s_i}(n',u)$ or is disjoint of it, and so that
\begin{equation}
\|x^{(3)}\|_\ell\!=\!\frac1{f(\ell)}\sum_{j=1}^\ell \|F_j(x^{(3)})\|\!
+\!\ell\vp(n')\!\max_{0\le i\le \!N+n+1} |\alpha_i|\!\le\!
 \frac1{f(\ell)}\!\sum_{j=1}^\ell \|F_j(x^{(3)})\|\!+2^{-n-i_0-1}\|x\|.
\end{equation}
To each of the  $F_j(x^{(3)})$ we can apply the induction hypothesis and 
obtain a splitting of each 
$$ \tilde F_j= \bigcup\{\supp(y_{t_i}(n',u)): i\ge i_0+3, u\le p(n')
 \text{ and } \supp(y_{s_i}(n',u))\subset F_j \} $$ 
($\tilde F_j$ is the ``$\tilde x^{(3)}$ version of $x^{(3)}$'')
into $\tilde G_j$ and $\tilde H_j$, with $\tilde G_j<\tilde H_j$
so that
\begin{align} \label{E:3.2.11}
\|x^{(3)}\|_\ell&\le 2^{-n-i_0-1} \|x\|+  
 \frac1{f(\ell)}c(n')
\sum_{j=1}^\ell \nl\tilde G_j(\tilde x^{(3)})\nr+ \nl\tilde H_j(\tilde x^{(3)})\nr\\
&\le2^{- n-i_0-1} \|x\|+ c(n')\nl \tilde x^{(3)}\nr\notag
\end{align}
(for the last inequality recall the definition of $\nl\cdot\nr$ and the fact that $g(2\ell)=f(\ell)$).

Finally, putting (\ref{E:3.2.7}), (\ref{E:3.2.9}) and (\ref{E:3.2.11}) together, observing that 
  $\|\cdot\|\le\nl\cdot\nr$,  the fact that $c(n)\ge 2^{-n-i_0}+ 2^{-n-i_0-1}+\vp(n+1)+ c(n')$ we obtain
\begin{align*}
\|x\|=\|x\|_\ell&\le \|x^{(1)}\|_\ell+\|x^{(2)}\|+\|x^{(3)}\|_\ell\\
    &\le \|\tilde x\| 2^{-n-i_0} + (1+\vp(n+i_0))\|\tilde x^{(2)}\| 
      + 2^{- n-i_0-1}\|\tilde x\| +c(n')\nl \tilde x^{(3)}\nr\\
     &c(n')\sup_{E<F} (\nl E(\tilde x)\nr+\nl F(\tilde x)) + [2^{-n-i_0}+ 2^{-n-i_0-1}+\vp(n+1)]\|\tilde x\|\\
   &\le     c(n)\sup_{E<F} (\nl E(\tilde x)\nr+\nl F(\tilde x)\nr),
\end{align*}
which finishes the induction step and the proof of the theorem.
\end{proof}

\noindent{\em Proof of Theorem }\ref{T:3.1}.
 Let $(x_n)$ be a seminormailzed basic sequence in $S$
 whose closed linear span is complemented in $S$. Since
$S$ is reflexive we can assume, by passing to an appropriate
 subsequence, that for some $a\in S$, $x_n=a+z_n$ where
 $z_n$ is weakly null. Since $a$ is element of the closed
linear span of $(x_n)$ and  $(x_n)$ is a seminormailzed basic sequence
 it follows that $a=0$ and, thus that $(x_n)$ is semi normalized
and weakly null. By applying the usual perturbation argument
 it  we can assume that $(x_n)$ is a seminormalized
block sequence and therefore apply Theorem \ref{T:3.2}
 with $k_n=2n+1$ for $n\in\N$ to obtain a subsequence
 $(y_n)$ so that $(y_n)$, $(y_{2n+1})$ and $(y_{2n})$ are
equivalent. In particular it follows for $Y=[y_n:n\in\N]$ that
$Y\approx Y\oplus Y$ (the complemented sum of $Y$ with itself).
Since $S$ has a subsymmetric basis, it also follows that
$S\approx S\oplus S$. This means that we are in the position
to apply Pe{\l}czy\'nski's decomposition method \cite{P}, which
 is so elegant that we cannot restrain ourselves to repeat it here.
We write $S=U\oplus Y$ (note that with $[x_i:i\in\N]$ also
 $Y$ is complemented in $S$) and $Y=V\oplus S$.
Then it follows that 
$$S\approx U\oplus Y \approx U\oplus Y\oplus Y
\approx S\oplus Y\approx S\oplus V\oplus S\approx S\oplus V\approx Y,$$
which proves our claim.\hfill$\square$

{\footnotesize
\noindent
Department of Mathematics, University of South Carolina, Columbia, SC 29208,
giorgis@math.sc.edu

\noindent
Department of Mathematics, Texas A$\&$M University, 
College Station, TX 77843, schlump@math.tamu.edu
}

\end{document}